\documentclass[11pt]{article}
\usepackage{amssymb} \usepackage{amsmath}
\usepackage{xspace}
\usepackage[top= 1.1 in,left=1.1 in, right=1.1 in, bottom=1.1 in]{geometry}
\usepackage{color}
\usepackage{ulem}

\newcommand{\blackboard}[1]{\mathbb#1}

\newcommand{\bbR}{\blackboard{R}}

\newcommand\LB[1]{\label{#1}} 
\newcommand\BE[2]{\begin{#1} #2 \end{#1}}

\newcommand\EQ[2]{\BE{equation}{\LB{#1} #2}}

 \newcommand\EQN[2]{\BE{equation*}{\LB{#1} #2}}
 \newcommand\EQn[1]{\BE{equation*}{ #1}}

\newcommand \lm{\lambda}
\newcommand \Lm{\Lambda}

\newcommand\dl{\delta}

\newcommand\lmn{\lambda_n}

\newcommand\iy{\infty}

\newcommand{\Arg}{\operatorname{Arg}}

\newcommand{\cE}{\mathcal{E}}
\newcommand{\cH}{\mathcal{H}}

\def\l{\left}  \def\r{\right}
\def\a{\alpha}  \def\b{\beta}

\begin{document}

\title { Density of zeros of the Cartwright class functions and the
Helson--Szeg\"o type condition}
\author{
Sergei\,A.\, Avdonin $^{a,b}$,
 Sergei A. Ivanov$^{c}$\footnote{\small
  Corresponding author. E-mail: sergei.a.ivanov@mail.ru.}
  \\
 $^a${\it Department of Mathematics and Statistics,
  University of Alaska Fairbanks, }\\
 {\it Faibanks, AK 99755-6660, USA}\\
$^b$ {\it Moscow Center for Fundamental and Applied Mathematics,  Moscow 119991, Russia}
\\
 $^c${\it Marine Geomagnetic Investigation Laboratory,}\\
    {\it
St. Petersburg branch of the Pushkov Institute of
 Russian Academy of Sciences,}\\
  {\it University Embankment, 5 B,  199034, St. Petersburg, Russia  }
\\
}
\date{}
\maketitle

\

\noindent { \emph{Key words and phrases}:} Helson--Szeg\"o condition, upper uniform density, exponential Riesz bases.

\

\noindent 2010 \emph{Mathematics Subject Classification}: Primary 30H35; Secondary 30J10,  42C99.

\begin{abstract}
	 B.\,Ya.\,Levin has proved that zero set of a sine type function can be presented as a union of a finite number of separated sets, that is an important result in the theory of exponential Riesz bases. In the present paper we extend Levin's result to a  more general class of entire functions $F(z)$  with
zeros in a strip  $0<q \leq \Im z \leq Q,$  such that  $|F(x)|^2$  satisfies the Helson--Szeg\"o condition. Moreover, we demonstrate that instead of the last condition one can require that $\log|F(x)|$  belongs to the BMO class.
\end{abstract}

\section{Introduction}

The theory of Riesz bases of complex exponentials
$\cE=\{e^{i\lmn t}\}$ or, in other words, nonharmonic Fourier series in $L^2$ on an interval began with the classical work by Paley and Wiener
\cite{PW}, which has motivated a great
deal of work by many mathematicians (see, e.g. the  references in \cite{AI95,KNP,Nik,Young}).  In \cite{L1} (see also \cite{L2}, Lectures 21--23)  Levin
developed the set of techniques that allows to connect the basis property of the family $\cE$ with the properties of an entire function with zero set $\{\lmn\}$.
Following this approach and attracting geometry of the Hardy space $H^2_+$ in the upper halfplane, Pavlov \cite{P} obtained the full description of the exponential Riesz bases.  For various generalizations of Pavlov's result see
\cite{N}, \cite{KNP},  \cite{M}, and  \cite{LS}.

In this paper we denote by $\Lm=\{\lm_n\}$  a sequence in $\mathbb{C}$ (multiple points are allowed) lying in a strip $S_{\Lm}$ parallel to the real axis:  $\sup |\Im \lm_n| < \iy$ . Without loss of generality of our results (and for convenience of notations) we assume also that this strip is situated in the upper half-plane, i.e. $\inf \Im \lm_n>0.$ The sequence
$\Lm$ is called separated, or uniformly discrete, if $\inf_{k \neq n} |\lm_k - \lm_n| >0.$
We say that $\Lm$ is relatively uniformly discrete if it can be decomposed
into a finite number of uniformly discrete subsequences.

The following notion plays an important role in the theory of  exponential bases. A function $F$ of the exponential type is called the sine type function if its zeros  $\lmn$ lie in a strip $S_{\Lm}$ and $|F(x)| \asymp 1,\, x \in \mathbb{R}.$ Levin \cite{L1} and Golovin \cite{G} proved the following
\BE{proposition}{\LB{LG} The family $\cE = \{e^{i\lmn t}\}$ forms a Riesz basis in $L^2(0,T)$ if
there is a sine type function with indicator diagram of width $T$ and uniformly discrete zero set $\Lm.$	}

In the case of a sine type function $f$ the  width  of the indicator diagram is the sum of exponential types $f$ in the upper and the lower halfplanes.

Pavlov's result can be formulated as the following statement.
\BE{proposition}{\LB{P}  The family $\cE$ forms a Riesz basis in $L^2(0,T)$ if and only if
		there is an entire function $F$ of exponential type
 with indicator diagram of width $T$ and uniformly discrete zero set $\Lm$	such that $|F(x)|^2$ satisfies the Helson--Szeg\"o condition: functions
 $
 u, v \in L^{\infty}(\mathbb{R}),
 {\Vert v \Vert}_{L^{\infty}(\mathbb{R})} < \pi/2,
 $
 can be found such that
 \EQN{HS}{
 	|F(x)|^2=\exp\{u(x)+\cH v(x)\}.
 }}
Here the map $v \mapsto \cH v$
denotes the Hilbert transform
for bounded functions:
$$
\cH v(x)=\frac 1\pi \, p.v.
\int_{-\infty}^\infty v(t)
\left \{ \frac 1{x-t}+\frac{t}{t^2+1} \right \} dt.
$$

The function $F$ in this theorem is called the generating function of the family
$\cE=\left\{e^{i\lambda_nt}\right\}$ on the interval $(0,T).$
This notion
plays a central role in the modern theory of
nonharmonic
Fourier series \cite{KNP,AI95}.	

If the sequence $\Lm$ is not uniformly discrete but relatively uniformly discrete, the family
$\cE$ is not a  Riesz basis of exponentials but it may
 form a Riesz basis of finite-dimensional subspaces (basis with brackets). A first result of such a type was obtained by Levin \cite{L1}. He proved that zero set $\Lm$ of a sine type function $F$ with indicator diagram of width $T$ is relatively uniformly discrete, and the corresponding family $\cE$ forms a Riesz basis with brackets in $L^2(0,T).$ Necessary and sufficient conditions for the Riesz basis property with brackets were formulated as  the Helson--Szeg\"{o}\ condition for the generating function $F$ provided that its zero set
$\Lm$ is relatively uniformly discrete \cite{AI01}. In the present paper we demonstrate that the last condition is excessive: if $F$ satisfies the Helson--Szeg\"{o}\ condition then its zero set is relatively uniformly discrete. Moreover, we prove that the conclusion is true if $\log|F(x)|$  belongs to the BMO class. This condition is weaker than  the Helson--Szeg\"{o}\ condition for $F$ (see the details below).

From Riesz bases of subspaces we can return to bases of individual functions if
 instead of exponentials consider the family of exponential divided differences. The theory of Riesz bases of exponential divided differences was developed in
\cite{AI01}, \cite{AM}. This theory has important applications to control theory of hybrid systems and delayed equations of neutral type, see, e.g. \cite{AM1}, \cite{IV}. In applications of Riesz bases of the exponential divided differences one needs to check two conditions:

(i) the generation function of the family $\cE$ with zero set $\Lm$ satisfies the Helson--Szeg\"o condition;

(ii)  $\Lm$ is relatively uniformly discrete. \\
The proof of the statement (ii) is typically rather complicated, since it requires a detailed information about the asymptotic behavior of the sequence $\lmn$. In the present paper we prove that the statement (i) implies (ii).

We need to introduce several more definitions. The upper uniform density $D_+(\Lm)$ of the sequence $\Lm$ is defined by
the formula
 \EQN{D}{
 	D_+(\Lm)=\lim_{r\to+\iy}\sup_{x\in \bbR}
 	\frac{\#\{\Re \Lm \cap [x,x+r)\}}r\,.
 }
This notion, as well as a similar notion of the lower uniform density $D_-(\Lm),$
 plays an important role in the theory of exponential bases $\cE$
 on an interval, see \cite{Seip}, \cite{AI01}, \cite{AM}.
If      $D_+(\Lm)=C<\iy$,  then in any strip $I \times i \,\mathbb{R}$
there are at most $C\,|I|$ points from $\Lm,$
that is equivalent to the fact that  $\Lm$ is relatively uniformly discrete.
Here $I$ is an interval of the real axis with length $|I|.$ 

\BE{definition}{   An entire function $F$ of exponential
type belongs to the Cartwright class  if
$$
     \int_{\bbR}  \frac{ \max \{ \log |F(x)|, 0 \} }{1+x^2} \, dx < \infty.
$$
}
The Cartwright class functions are studied in detail in the classical monograph \cite{Levin}.

\BE{definition}{ The mean oscillation $p_I\,(f)$ of a locally integrable function $f$ over an interval $I$ is defined as
\begin{equation*} \label{pfi}
p_I\,(f)= {\frac {1}{|I|}}\int _{I} \left|f(x)- f_I \, \right|\,  dx\,,
\end{equation*}
where $f_I$ is the mean value of $f$ on $I:$
$$
f_I=\frac1{|I|}\int_I f\,dx.
$$
 A locally integrable function on $\mathbb{R}$  belongs to the BMO class if the supremum
of its mean oscillation $\sup_{I\in {\cal J}}\,p_I\,(f),$
taken over the set of all intervals of the real axis, is finite.}

It is known that a locally integrable function $f$ on $\bbR$ belongs to the BMO class if and only if it can be written as
\EQn{
  f=u+\cH v,
  }
where $u,v \in L^\iy(\mathbb{R})$.

In the present paper we demonstrate that if
the function $F$ with zero set  $\Lm$   belongs to the Cartwright class   and
$$\log|F(x)|\in  \text{BMO},$$
then $\Lm$ has a finite upper uniform density:
$D_+(\Lm)<\iy$.

It is convenient to prove this result in a contrapositive form:
\BE{theorem}{\LB{main} Let the function $F$  with zero set
	$\Lm$ in a strip $ 0 < \inf \Im z \leq \sup \Im z<\iy$ belong to  the Cartwright class and
$D_+(\Lm)=\iy.$
Then
\EQN {notBMO}{
\log|F(x)|\notin  \text{BMO}.
}
}
In particular, it implies that $F$ does not satisfy the Helson--Szeg\"o condition.

\BE{remark} {The authors are grateful to the anonymous referee for examples of the Cartwright class functions  with zero set of infinite upper uniform density
	lying in a strip parallel to the real axis. 

The first example demonstrates a function with zeros of unbounded multiplicity:	
	$$
	f(z)=\prod_{n=1}^\iy \cos^n(z/n^3).
	$$
	
One may get rid of the multiplicity condition.  For that, in the second example we consider the function	
	$$
	f_n(z)=\cos\left [\frac\pi2(3^{-n}+3^{-n^2})z \right],
	$$ 	
	vanishing at the points 	
	$$
	z_{kn}=3^k-\frac{3^{k-n^2}}{3^{-n}+3^{-n^2}}
	$$
	for every integers $k,\, n,$ such that $k>n$. It is easy to check that the function
	$$
	f(z)=\prod_{n=1}^\iy f_n(z)
	$$
	has at least  $k/2$ simple zeros on the interval
	 $(3^k-1,3^k),$ for all sufficiently large $k.$
}

\section{Proof of Theorem \ref{main}}

Suppose that $D_+(\Lm)=\iy$.
Following \cite{log}, 
we introduce  a continuous branch $\varphi_z(t)$ of the $\Arg b_z(t)$ for the   Blaschke factor  
$$
b_z(t)=\frac{1-t/z}{1-t/\bar z}.
$$
For this purpose we
put 
$$
\psi_z(t)=\arctan \frac{yt}{|z|^2-xt}, \  z=x+iy, \  y>0,
$$
and set
\EQN{phi+}{
\varphi_z(t)=\BE{cases}
{\psi_z(t)+\pi, \  &t>\frac {|z|^2}x,\\
\pi/2,                                             &t=\frac {|z|^2}x,\\
\psi_z(t), \  &t<\frac {|z|^2}x;
}
}
if $x>0,$ and
\EQN{phi-}{
\varphi_z(t)=\BE{cases}
{\psi_z(t), \       &t>\frac {|z|^2}x,\\
-\pi/2,                                             &t=\frac {|z|^2}x,\\
\psi_z(t) -\pi, \     &t<\frac {|z|^2}x.
}
}
if $x<0.$
The following result has been obtained in \cite{log}.
\BE{proposition}{
\LB{logC}
Let the function $F$  with zero set $\{z_n\}$
 belong to the Cartwright class. Then the Hilbert transform $\cH$ of $\log|F|$ is presented in the following form:
\begin{equation} \label{hlf}
\cH \,[\,\log|F|\,](t)=\theta+\frac T2 \, t-\sum_n \varphi_{z_n}(t),
\end{equation}
where $T$ is the width  of the indicator diagram of $F$ and $\theta$ is a constant.
}
Note that the series converges because the sequence $\{z_n\}$ satisfies the Blaschke condition
$$
\sum \frac{\Im z_n}{|z_n|^2}<\iy.
$$

As the first step in proving Theorem \ref{main} we demonstrate  that  the sum
$$
\Phi(t)=\sum_n \varphi_{z_n}(t)
$$
rapidly increases if $D_+(\Lm)=\iy\,. $

\BE{lemma}{\LB{fast0} If $D_+(\Lm)=\iy$, then for any $M>0$ there exists $a$ such that
\EQN{fast1}{
\Phi(a+1)-\Phi(a) \geq M.
}
}
\emph{Proof} of the lemma. Direct calculation gives
$$
\psi'_{x+iy}\,(t)=\frac{1}{1+\l[\frac{yt}{(x^2+y^2)-xt} \r]^2}\,\cdot\,
\l[ \frac{y}{(x^2+y^2)-xt}+\frac{yxt}{\l[\frac{yt}{(x^2+y^2)-xt} \r]^2 } \r]
$$
$$
=\frac{y(x^2+y^2)}{(x^2+y^2-xt)^2+y^2t^2}.
$$

For  $t=x+\dl,\,$  $|\dl| \leq 1,$ we have
$$
\psi'_{x+iy}\,(t)=\frac y{y^2+\dl^2}.
$$
For  $z$ is in the strip $0<\a \le \Im z \le \b<\iy$ the minimum of the function $ y/(y^2+\dl^2)$ is positive:
\EQ{const}{
\psi'_{x+iy}\,(t)\ge \frac y{y^2+1}\ge \min\l\{ \frac\a{\a^2+1}, \frac \b{\b^2+1}\r\} =c>0 \quad {\rm{for}} \ t \in [x-1,x+1].
}

Now we put
$$\Lm_a = \{z_n|\, \Re z_n \in [a,a+1]  \}
$$
and choose a segment $[a, a+1]$ such that $\# \Lm_a \geq M/c ;$ it is possible by the assumption
of the lemma.
Evidently,
$$
\Phi(a+1)-\Phi(a) \geq
\sum_{z_n\in \Lm_a} \varphi_{z_n}(a+1) -\sum_{z_n\in \Lm_a} \varphi_{z_n}(a)=
\sum_{z_n\in \Lm_a}[ \varphi_{z_n}(a+1) - \varphi_{z_n}(a)].
$$
For every function $\varphi_{z_n}$ in the RHS of this inequality, \eqref{const} gives
$$
\varphi_{z_n}(a+1) - \varphi_{z_n}(a)\ge c,
$$
and therefore,
$$
\Phi(a+1)-\Phi(a)\ge   M.
$$
The lemma is proved.

\BE{lemma}{\LB{fast2} Let $g$ be an increasing function on $\bbR$ and,
for any $M>0,$ there exists $a \in \mathbb{R}$ such that
\EQ{fast3}{
g(a+1)-g(a) \geq M.
}
Then
$$
\sup_{I=[a-1,a+2]} p_I\, (f)=\iy,
$$
and therefore, $g$ does not belong to BMO
}
\emph{Proof} of this lemma  is  based on  straightforward calculations.
Let $M$ and $a$ satisfy the inequality \eqref{fast3}.
Set $I=I_M=[a-1,a+2]$. We will demonstrate  that
\EQ1{
p_I(g)\ge\frac M6.
}
Evidently, the mean value $g_I$ belongs to the interval $[g(a-1), g(a+2)]$ and
\EQN2{
g(a-1)\le g(a)\le  g(a)+M\le  g(a+1)\le  g(a+2).
}
 We consider two cases:

(i) $g_I\in [g(a-1), g(a+\frac M2))$, \ \ (ii) $g_I\in [g(a+\frac M2), g(a+2)]$.

\textbf{Case (i).} \ If $x\in [a+1,a+2]$, then
$$
g(x) \geq g(a+1)\ge g(a)+M \geq  g(a)+\frac M2 \geq g_I
$$
and
$$
|g(x)-g_I|= g(x) - g_I \ge  g(a+1) - g_I  \ge g(a)+M-g_I\ge \frac M2.
$$
Therefore,
$$
p_I(g)\ge \frac13\int\limits_{a+1}^{a+2}|g(x)-g_I|dx\ge \frac M6.
$$

\textbf{Case (ii).} \  If $x\in [a-1,a]$, then
$$
g_I \geq g(a)\ge g(x)
$$
and
$$
|g(x)-g_I|=-g(x)+g_I\ge g(a)+\frac M2 -g(x)\ge g(a)+\frac M2 -g(a)\ge \frac M2.
$$
Therefore,
$$
p_I(g)\ge \frac13\int\limits_{a-1}^{a}|g(x)-g_I|dx\ge \frac M6.
$$
We obtained (\ref1) and the lemma is  proved.

Now we are able to complete the proof of Theorem \ref{main}.
The functional $p_I$ of the term $Tt/2$ is constant  on intervals of length 3.
Then from (\ref{hlf}) and Lemmas \ref{fast0} and \ref{fast2} it follows that
\EQN{ImlgBMO}{
	\cH[\, \log|F|\,]\not \in \text{BMO}.
}
The Hilbert transform maps the  \text{BMO} class to itself and $\cH^2=-Id$.
So,
$$
\log|F|=-\cH\cH \, [\,\log|F|\,]
$$
and therefore, $\log|F|$ does not belong to the BMO. The theorem is proved.

\BE{remark}{
In fact, to prove the main theorem, the explicit  formulas  in Proposition \ref{logC} are not necessary.
 In the case of infinite density we only need  to show the ``rapid growth"  of the continuous branch of Arg F (in the sense Lemma \ref{fast0}).
 It has been proved by the authors  in other way, but the result of the work \cite{log}  significantly simplifies the proof.

Let us also note that the connections between the  continues branches of the argument of the Blaschke product, the BMO space, and the basis property of exponential families were discussed in  \cite[Section I.4]{KNP}.
}
\vskip3mm

\noindent {\bf  Acknowledgments}\\
The research of Sergei Avdonin was  supported  in part by the National Science Foundation,
grant DMS 1909869.

\end{document}